\documentclass[12pt,a4paper]{article}
\usepackage{amssymb}

\usepackage{amsmath}

\newtheorem{theorem}{Theorem}[section]

\newtheorem{lemma}[theorem]{Lemma}

\newtheorem{example}[theorem]{Example}

\begin{document}
\title{Composition-Diamond Lemma for Tensor Product of Free Algebras\footnote{Supported by the
NNSF of China (No.10771077) and the NSF of Guangdong Province
(No.06025062).}}

\author{
L. A. Bokut\footnote {Supported by the RFBR  and
the Integration Grant of the SB RAS (No. 1.9).} \\
{\small \ School of Mathematical Sciences, South China Normal
University}\\
{\small Guangzhou 510631, P. R. China}\\
{\small Sobolev Institute of Mathematics, Russian Academy of
Sciences}\\
{\small Siberian Branch, Novosibirsk 630090, Russia}\\
{\small Email: bokut@math.nsc.ru}\\
\\
 Yuqun
Chen\footnote {Corresponding author.} \  and Yongshan Chen\\
{\small \ School of Mathematical Sciences, South China Normal
University}\\
{\small Guangzhou 510631, P. R. China}\\
{\small Email: yqchen@scnu.edu.cn}\\
{\small chenyongshan@gmail.com}}

\date{}
\maketitle

\maketitle \noindent\textbf{Abstract:} In this paper, we establish
Composition-Diamond lemma for tensor product $k\langle X\rangle
\otimes k\langle Y\rangle$ of two free algebras over a field. As an
application, we construct a Gr\"{o}bner-Shirshov basis in $k\langle
X\rangle \otimes k\langle Y\rangle$ by lifting a
Gr\"{o}bner-Shirshov basis in $k[X] \otimes k\langle Y\rangle$,
where $k[X]$ is a commutative algebra.


\noindent {\bf AMS} Mathematics Subject Classification(2000): 16S10,
16S15, 13P10

\section{Introduction}

In 1962, A. I. Shirshov \cite{s62} invented a theory of
one-relator Lie algebras $Lie(X|s=0)$ that was in a full analogy,
by statements, but not the method, of celebrated  Magnus's theory
of one-relater groups \cite{mag30} and \cite{mag32}, see also
\cite{mksbook} and \cite{lsbook}. In particular, A. I. Shirshov
proved the algorithmic decidality of the word problem for any
one-relator Lie algebra. To do it, he created a theory that is now
called the Gr\"{o}bner-Shirshov bases theory for Lie algebras
$Lie(X|S)$ presented by generators and defining relations. The
main technical notion of the latter Shirshov's theory was a notion
of composition $(f,g)_w$ of two Lie polynomials, $f,g\in Lie(X)$
relative to some associative word $w$. Based on it, he defined an
infinite algorithm of adding to some set $S$ of Lie polynomials
all non-trivial compositions until one will get a set $S^*$ that
is closed under compositions, i.e., any non-trivial composition of
two polynomials from $S^*$ belongs to $S^*$ (and leading
associative words $\bar{s}$ of polynomials $s\in S^*$ do not
contain each others as subwords). In addition, $S$ and $S^*$
generated the same ideal, i.e., $Id(S)=Id(S^*)$. $S^*$ is now
called a Gr\"{o}bner-Shirshov basis of $Id(S)$. Then he proved the
following lemma.

\emph{Let $Lie(X)\subset k\langle X\rangle$  be a free Lie algebra
over a field $k$ viewed as the algebra of Lie polynomials in the
free algebra $k\langle X\rangle$, and $S$ a subset in $Lie(X)$. If
$f\in Id(S)$, then $\bar{f}=u\bar{s}v$, where $s\in S^*, \ u,v\in
X^*, \ \bar{f},\bar{s}$ are leading associative words of Lie
polynomials $f,s$ correspondingly, and $X^*$ the free monoid
generated by $X$.}

He used the following easy corollary of his lemma.

\emph{$Irr(S)=\{[u] \ | \  u\neq{a\bar{s}b}, \ s\in{S},\
a,b\in{X^*}\}$ is a linear basis of the algebra
$Lie(X|S)=Lie(X)/Id(S)$, where $u$ is an associative Lyndon-Shirshov
word in $X^*$ and $[u]$ the corresponding non-associative
Lyndon-Shirshov word under Lie brackets $[xy]=xy-yx$.}

To define the Lie composition $(f,g)_w$ of two, say, monic Lie
polynomials, where $\bar{f}=ac, \ \bar{g}=cb, \ c\neq 1, \ a,b,c$
are associative words, and $w=acb$, A. I. Shirshov defines first
the associative composition $fb-ag$. Then he puts on $fb$ and $ag$
special brackets $[fb],[ag]$ in according with his paper
 \cite{s58}. The result is $(f,g)_w=[fc]-[cg]$. Following
 \cite{s62}, one can easily get the same lemma for a free
associative algebra: Let $S\subset k\langle X\rangle$ and $S^*$ be
as before. If $f\in Id(S)$, then $\bar{f}=a\bar{s}b$ for some
$s\in S^*, \ a,b\in X^*$. It was formulated lately by L. A. Bokut
 \cite{b76} as an analogy of Shirshov's Lie composition lemma, and
by G. Bergman \cite{b} under the name ``Diamond lemma" after
celebrated Newman's Diamond Lemma for graphs \cite{newman}.

Shirshov's lemma is now called the Composition-Diamond lemma for Lie
and associative algebras. Its nowadays formulation see, for example,
in the next section in this paper.

Independently this kind of ideas were discovered by H. Hironaka
 \cite{H64} for the power series algebras and by B. Buchberger
 \cite{bu65}, \cite{bu70} for the polynomial algebras. B. Buchberger
suggested the name ``Gr\"{o}bner bases". It is well known and well
recognized that applications of Gr\"{o}bner  bases in mathematics
(particulary, in algebraic geometry), computer science and
informatics are innumerable large.

At present, there are quite a few Compositions-Diamond lemmas
(CD-lemma for short) for different classes of non-commutative and
non-associative algebras. Let us mention some.

A. I. Shirshov \cite{s62'} proved himself CD-lemma for commutative
(anti-commutative) non-associative algebras, and mentioned that this
lemma is also valid for non-associative algebras. It gives solution
of the word problems for these classes of algebras. For
non-associative algebras, this (but not CD-lemma) was known, see A.
I. Zhukov \cite{Z}.

 A.~A.~Mikhalev \cite{MikhPetr} proved a CD-lemma for Lie super-algebras.

 T. Stokes \cite{st90} proved a CD-lemma for left ideals of an algebra $k[X]\otimes E_k(Y)$,
 the tensor product of Exterier (Grassman) algebra and a polynomial algebra.

A. A. Mikhalev and E. A. Vasilieva \cite{MV} proved a CD-lemma for
the free supercommutative polynomial algebras.

 A. A. Mikhalev and A.~A. Zolotykh  \cite{MZ} proved a CD-lemma
 for $k[X]\otimes k\langle Y\rangle$, the tensor product of a polynomial algebra and a free algebra.

L. A. Bokut, Y. Fong and W. F. Ke \cite{bfk} proved a CD-lemma for
associative conformal algebras.

L. Hellstr\"{o}m \cite{H} proved a CD-lemma for a non-commutative
power series algebra.

 S.-J. Kang and K.-H. and Lee \cite{kl1}, \cite{kl3} and E. S. Chibrikov \cite{ch} proved a
CD-lemma for a module over an algebra.

D. R. Farkas, C. D. Feustel and E. L. Green \cite{ffg93} proved a
CD-lemma for path algebras.

L. A. Bokut and K. P. Shum \cite{b06} proved a CD-lemma for
$\Gamma$-algebras.

Y. Kobayashi \cite{kobayashi}  proved a CD-lemma for algebras
based on well-ordered semigroups, and L. A. Bokut, Yuqun Chen and
Cihua Liu \cite{bcl} proved a CD-lemma for dialgebras (see also
\cite{bc}).

Let $X$ and $Y$ be sets and $k\langle X\rangle \otimes k\langle
Y\rangle$ the tensor product algebra. In this paper, we give the
Composition-Diamond lemma for the algebra $k\langle X\rangle \otimes
k\langle Y\rangle$.  Also we will prove a theorem on the pair of
algebras $(k[X]\otimes k\langle Y\rangle, k\langle X\rangle\otimes
k\langle Y\rangle)$ in spirit of Eisenbud, Peeva and Sturmfels
theorem \cite{eisenbud} on $(k[X],k\langle X\rangle)$.

\section{Preliminaries}

We cite some concepts and results from the literature (\cite{s62},
\cite{b72}, \cite{b76}) concerning with the Gr\"{o}bner-Shirshov
bases theory of associative algebras.

Let $k$ be a field, $k\langle X\rangle$ the free associative algebra
over $k$ generated by $X$ and $ X^{*}$ the free monoid generated by
$X$, where the empty word is the identity which is denoted by 1. For
a word $w\in X^*$, we denote the length of $w$ by $|w|$.

A well order $>$ on $X^*$ is monomial if it is compatible with the
multiplication of words, that is, for $u, v\in X^*$, we have
$$
u > v \Rightarrow w_{1}uw_{2} > w_{1}vw_{2},  \ for \  all \
 w_{1}, \ w_{2}\in  X^*.
$$
A standard example of monomial order on $X^*$ is the deg-lex order
to compare two words first by degree and then lexicographically,
where $X$ is a linearly ordered set.

Let $f\in k\langle X\rangle$ with the leading word $\bar{f}$. We say
that $f$ is monic if $\bar{f}$ has coefficient 1.

Let $f$ and $g$ be two monic polynomials in \textmd{k}$\langle
X\rangle$ and $<$ a well order on $X^*$. Then, there are two kinds
of compositions:

$(1)$ If \ $w$ is a word such that $w=\bar{f}b=a\bar{g}$ for some
$a,b\in X^*$ with $|\bar{f}|+|\bar{g}|>|w|$, then the polynomial
 $(f,g)_w=fb-ag$ is called the intersection composition of $f$ and
$g$ with respect to $w$.

$(2)$ If  $w=\bar{f}=a\bar{g}b$ for some $a,b\in X^*$, then the
polynomial $(f,g)_w=f - agb$ is called the inclusion composition of
$f$ and $g$ with respect to $w$.

Let $S\subset$ $\textmd{k}\langle X\rangle$ with each $s\in S$
monic. Then the composition $(f,g)_w$ is called trivial modulo
$(S,w)$ if $(f,g)_w=\sum\alpha_i a_i s_i b_i$, where each
$\alpha_i\in k, \ a_i,b_i\in X^{*}, \ s_i\in S$ and $\overline{a_i
s_i b_i}<w$. If this is the case, then we write
$$
(f,g)_w\equiv0\quad mod(S,w).
$$
In general, for $p,q\in k\langle X\rangle$, we write $p\equiv q\quad
mod(S,w)$ which means that $p-q\equiv0\quad mod(S,w)$.

We call the set $S$ with respect to the monomial order $<$ a
Gr\"{o}bner-Shirshov basis in $k\langle X\rangle$ if any composition
of polynomials in $S$ is trivial modulo $S$.

\begin{lemma}\label{l1.4}
(Composition-Diamond lemma for associative algebras) \ Let $S\subset
k \langle X\rangle$ be a set of monic polynomials and $<$ a monomial
order on $X^*$. Then the following statements are equivalent:
\begin{enumerate}
\item[(1)] $S $ is a Gr\"{o}bner-Shirshov basis in $k \langle X\rangle$.
\item[(2)] $f\in Id(S)\Rightarrow \bar{f}=a\bar{s}b$
for some $s\in S$ and $a,b\in  X^*$, where $Id(S)$ is the ideal of
$k \langle X\rangle$ generated by $S$.
\item[(3)] $Irr(S) = \{ u \in X^* |  u \neq a\bar{s}b ,s\in S,a ,b \in X^*\}$
is a basis of the algebra $A=k\langle X | S \rangle$.
\end{enumerate}
\end{lemma}

\section{Composition-Diamond Lemma for Tensor Product}

Let $X$ and $Y$ be linearly ordered sets, $T=\{yx=xy|x\in X, \ y\in
Y\}$. With the deg-lex order ($y>x$ for any $x\in X, \ y\in Y$) on
$(X\cup Y)^*$, $T$ is a Gr\"{o}bner-Shirshov basis in $k\langle
X\cup Y\rangle$. Then, by Lemma \ref{l1.4}, the set
$$
N=X^*Y^*=Irr(T)=\{u=u^Xu^Y| u^X\in X^* \ and \ u^Y\in Y^*\}
$$
is the normal words of the tensor product
$$
k\langle X\rangle \otimes k\langle Y\rangle=k\langle X\cup Y \ | \ T
\rangle.
$$

Let $kN$ be a $k$-space spanned by $N$. For any
$u=u^Xu^Y,v=v^Xv^Y\in N$, we define the multiplication of the normal
words as follows
$$
uv=u^Xv^Xu^Yv^Y\in N.
$$
Then,  $kN$ is  exactly tensor product algebra $k\langle X\rangle
\otimes k\langle Y\rangle$, that is, $kN=k\langle X\cup
Y|T\rangle=k\langle X\rangle \otimes k\langle Y\rangle$.

Let $``>_X"$ and $``>_Y"$ be any monomial orders on $X^*$ and $Y^*$
respectively.  Now, we order the set $N$. For any
$u=u^Xu^Y,v=v^Xv^Y\in N$,
$$
u>v\Leftrightarrow u^X>_Xv^X \  or \ (u^X=v^X \ and \ u^Y>_Yv^Y).
$$
It is obvious that $>$ is a monomial order on $N$. Such an order is
also called the deg-lex order on $N=X^*Y^*$. We will use this order
in the sequel unless others stated.

For any polynomial $f\in k\langle X\rangle \otimes k\langle
Y\rangle$, $f$ has a unique presentation of the form
$$
f=\alpha_{\bar{f}}\bar{f}+\sum\alpha_iu_i,
$$
where $\bar{f},u_i\in N,\bar{f}>u_i,\alpha_{\bar{f}},\alpha_i\in k.$

The proof of the following lemma are straightforward.

\begin{lemma}\label{l2} Let $f\in k\langle X\rangle \otimes k\langle Y\rangle$ be a monic polynomial. Then
$\overline{ufv}=u\bar{f}v$ for any $u,v\in N$.
\end{lemma}

Now, we give the definition of compositions. Let $f$ and $g$ be
monic polynomials of $k\langle X\rangle \otimes k\langle Y\rangle$
and $w=w^Xw^Y\in N$. Then we have the following compositions.

$1.$  Inclusion

$1.1$ $X$-inclusion only

Suppose that $w^X=\bar{f}^X=a\bar{g}^Xb$ for some $a,b\in X^*$, and
$\bar{f}^Y, \ \bar{g}^Y$ are disjoint. Then there are two
compositions according to $w^Y=\bar{f}^Yc\bar{g}^Y$ and
$w^Y=\bar{g}^Yc\bar{f}^Y$ for  $c\in Y^*$, respectively:
$$
(f,g)_{w_1}=fc\bar{g}^Y-\bar{f}^Ycagb, \ \
w_1=f^X\bar{f}^Yc\bar{g}^Y
$$
and
$$
(f,g)_{w_2}=\bar{g}^Ycf-agbc\bar{f}^Y, \ \
w_2=f^X\bar{g}^Yc\bar{f}^Y.
$$

$1.2$  $Y$-inclusion only

Suppose that $w^Y=\bar{f}^Y=c\bar{g}^Yd$ for  $c,d\in Y^*$, and
$\bar{f}^X, \ \bar{g}^X$ are disjoint.  Then there are two
compositions according to $w^X=\bar{f}^Xa\bar{g}^X$ and
$w^X=\bar{g}^Xa\bar{f}^X$ for  $a\in X^*$,  respectively:
$$
(f,g)_{w_1}=fa\bar{g}^X-\bar{f}^Xacgd, \ \
w_1=\bar{f}^Xa\bar{g}^Xf^Y
$$
and
$$
(f,g)_{w_2}=\bar{g}^Xaf-cgda\bar{f}^X \ \
w_2=\bar{g}^Xa\bar{f}^Xf^Y.
$$

 $1.3$ $X,Y$-inclusion

 Suppose that $w^X=\bar{f}^X=a\bar{g}^Xb$ for some $a,b\in X^*$ and
$w^Y=\bar{f}^Y=c\bar{g}^Yd$ for some $c,d\in Y^*$. Then
$$(f,g)_w=f-acgbd.$$

The transformation $f\mapsto (f,g)_w=f-acgbd$ is called the
\emph{elimination of the leading word}
(ELW) of $g$ in $f$.

$1.4$ $X,Y$-skew-inclusion

Suppose that $w^X=\bar{f}^X=a\bar{g}^Xb$ for some $a,b\in X^*$ and
$w^Y=\bar{g}^Y=c\bar{f}^Yd$ for some $c,d\in Y^*$. Then
$$(f,g)_w=cfd-agb.$$

$2.$  Intersection

$2.1$ $X$-intersection only

Suppose that $w^X=\bar{f}^Xa=b\bar{g}^X$ for some $a,b\in X^*$ with
$|\bar{f}^X|+|\bar{g}^X|>|w^X|$, and $\bar{f}^Y, \ \bar{g}^Y$ are
 disjoint.  Then there
 are two compositions  according to $w^Y=\bar{f}^Yc\bar{g}^Y$ and
$w^Y=\bar{g}^Yc\bar{f}^Y$ for  $c\in Y^*$, respectively:
$$
(f,g)_{w_1}=fac\bar{g}^Y-\bar{f}^Ycbg, \ \
w_1=w^X\bar{f}^Yc\bar{g}^Y
$$
and
$$
(f,g)_{w_2}=\bar{g}^Ycfa-bgc\bar{f}^Y, \ \
w_2=w^X\bar{g}^Yc\bar{f}^Y.
$$

$2.2$  $Y$-intersection only

Suppose that $w^Y=\bar{f}^Yc=d\bar{g}^X$ for some $c,d\in Y^*$ with
$|\bar{f}^Y|+|\bar{g}^Y|>|w^Y|$, and  $\bar{f}^X, \ \bar{g}^X$ are
disjoint. Then there
 are two compositions according to $w^X=\bar{f}^Xa\bar{g}^X$ and
$w^X=\bar{g}^Xa\bar{f}^X$ for  $a\in X^*$,  respectively:
$$
(f,g)_{w_1}=fca\bar{g}^X-\bar{f}^Xadg,\ \ w_1=\bar{f}^Xa\bar{g}^Xw^Y
$$
and
$$
(f,g)_{w_2}=\bar{g}^Xafc-dga\bar{f}^X, \ \
w_2=\bar{g}^Xa\bar{f}^Xw^Y.
$$

$2.3$ $X,Y$-intersection

If $w^X=\bar{f}^Xa=b\bar{g}^X$ for some $a,b\in X^*$ and
$w^Y=\bar{f}^Yc=d\bar{g}^Y$ for some $c,d\in Y^*$ together with
$|\bar{f}^X|+|\bar{g}^X|>|w^X|$ and $|\bar{f}^Y|+|\bar{g}^Y|>|w^Y|$,
then
$$(f,g)_w=fac-bdg.$$

$2.4$ $X,Y$-skew-intersection

If $w^X=\bar{f}^Xa=b\bar{g}^X$ for some $a,b\in X^*$ and
$w^Y=c\bar{f}^Y=\bar{g}^Yd$ for some $c,d\in Y^*$ together with
$|\bar{f}^X|+|\bar{g}^X|>|w^X|$ and $|\bar{f}^Y|+|\bar{g}^Y|>|w^Y|$,
then
$$(f,g)_w=cfa-bgd.$$

$3.$  Both inclusion and intersection

$3.1$ $X$-inclusion and $Y$-intersection

There are two cases to consider.

 If
$w^X=\bar{f}^X=a\bar{g}^Xb$ for some $a,b\in X^*$ and
$w^Y=\bar{f}^Yc=d\bar{g}^Y$ for some $c,d\in Y^*$  with
$|\bar{f}^Y|+|\bar{g}^Y|>|w^Y|$, then
$$(f,g)_w=fc-adgb.$$

If $w^X=\bar{f}^X=a\bar{g}^Xb$ for some $a,b\in X^*$ and
$w^Y=c\bar{f}^Y=\bar{g}^Yd$ for some $c,d\in Y^*$  with
$|\bar{f}^Y|+|\bar{g}^Y|>|w^Y|$, then
$$(f,g)_w=cf-agbd.$$

$3.2$  $X$-intersection and $Y$-inclusion

There are two cases to consider.

 If $w^X=\bar{f}^Xa=b\bar{g}^X$ for
some $a,b\in X^*$ with $|\bar{f}^X|+|\bar{g}^X|>|w^X|$ and
$w^Y=\bar{f}^Y=c\bar{g}^Yd$ for some $c,d\in Y^*$, then
$$(f,g)_w=fa-bcgd.$$

If $w^X=\bar{f}^Xa=b\bar{g}^X$ for some $a,b\in X^*$ with
$|\bar{f}^X|+|\bar{g}^X|>|w^X|$ and $w^Y=c\bar{f}^Yd=\bar{g}^Y$ for
some $c,d\in Y^*$, then
$$(f,g)_w=cfad-bg.
$$

From Lemma \ref{l2}, it follows that for any case of compositions
$$\overline{(f,g)_w}<w.$$

If $Y=\emptyset$, then the compositions of $f,g$ are the same in
$k\langle X\rangle$.

Let $S$ be a monic subset of $k\langle X\rangle \otimes k\langle
Y\rangle$ and $f,g\in S$. A composition $(f,g)_w$ is said to be
\emph{trivial modulo} $(S,w)$, denoted by
$$
(f,g)_w\equiv 0 \ \ mod(S,w), \ \mbox{ if } \
(f,g)_w=\sum_i\alpha_ia_is_ib_i,
$$
where $a_i,b_i\in N, \ s_i\in S,
  \ \alpha_i\in k$ and $a_i\bar{s_i}b_i<w$ for any  $i$.

Generally, for any $p,q\in k\langle X\rangle \otimes k\langle
Y\rangle,\ p\equiv q \ \ mod (S,w)$ if and only if $p-q\equiv 0 \ \
mod(S,w).$

 $S$ is called a \emph{Gr\"{o}bner-Shirshov basis} in
$k\langle X\rangle \otimes k\langle Y\rangle$ if all compositions of
elements in $S$ are trivial modulo $S$ and corresponding to $w$.

\begin{lemma}\label{l3}Let $S$ be a Gr\"{o}bner-Shirshov basis in
$k\langle X\rangle \otimes k\langle Y\rangle$ and $s_1,s_2\in S$. If
$w=a_1\bar{s_1}b_1=a_2\bar{s_2}b_2$ for some $a_i,b_i\in N, \
i=1,2$, then $a_1s_1b_1\equiv a_2s_2b_2 \ mod(S,w)$. \end{lemma}

\textbf{Proof: }There are four cases to consider.

\emph{Case 1} Inclusion

\emph{(1.1)} $X$-inclusion only

Suppose that $w_1^X=\bar {s_1}^X=a\bar{s_2}^Xb, \ a,b\in X^*$ and
$\bar {s_1}^Y, \ \bar{s_2}^Y$ are disjoint. Then $a_2^X=a_1^Xa$ and
$b_2^X=bb_1^X$. There are two cases to consider:
$w_1^Y=\bar{s_1}^Yc\bar{s_2}^Y$ and $w_1^Y=\bar{s_2}^Yc\bar{s_1}^Y$,
where $c\in Y^*$.

For $w_1^Y=\bar{s_1}^Yc\bar{s_2}^Y$, we have
$w_1=s_1^X\bar{s_1}^Yc\bar{s_2}^Y, \ a_2^Y=a_1^Y\bar{s_1}^Yc$,
$b_1^Y=c\bar{s_2}^Yb_2^Y$,
$w=a_1w_1b_1^Xb_2^Y=a_1\bar{s_1}ac\bar{s_2}b_2$ and
\begin{eqnarray*}
a_1s_1b_1-a_2s_2b_2&=&a_1s_1b_1^Xc\bar{s_2}^Yb_2^Y-a_1^Xaa_1^Y\bar{s_1}^Ycs_2bb_1^Xb_2^Y\\
&=&a_1(s_1c\bar{s_2}^Y-\bar{s_1}^Ycas_2b)b_1^Xb_2^Y\\
&=&a_1(s_1,s_2)_{w_1}b_1^Xb_2^Y\\
&\equiv&0 \ \ \ \ \ \ \ mod(S,w).
\end{eqnarray*}

For $w_1^Y=\bar{s_2}^Yc\bar{s_1}^Y$,
 we have
$w_1=s_1^X\bar{s_2}^Yc\bar{s_1}^Y, \ a_1^Y=a_2^Y\bar{s_2}^Yc$,
$b_2^Y=c\bar{s_1}^Yb_1^Y$, $w=a_1^Xa_2^Yw_1b_1$ and
\begin{eqnarray*}
a_1s_1b_1-a_2s_2b_2&=&a_1^Xa_2^Y\bar{s_2}^Ycs_1b_1-a_1^Xaa_2^Ys_2bb_1^Xc\bar{s_1}^Yb_1^Y\\
&=&a_1^Xa_2^Y(\bar{s_2}^Ycs_1-as_2bc\bar{s_1}^Y)b_1\\
&=&a_1^Xa_2^Y(s_1,s_2)_{w_1}b_1\\
&\equiv&0 \ \ \ \ \ \ \ mod(S,w).
\end{eqnarray*}

\emph{(1.2)} $Y$-inclusion only

This case is similar to (1.1).

\emph{(1.3)} $X,Y$-inclusion

We may assume that $\bar{s_2}$ is a subword of $\bar{s_1}$, i.e.,
$w_1=\bar{s_1}=ac\bar{s_2}bd$, $a,b\in X^*$, $c,d\in Y^*$,
 $a_2^X=a_1^Xa$,  $b_2^X=bb_1^X$, $a_2^Y=a_1^Yc$ and $b_2^Y=db_1^Y$. Thus,
 $a_2=a_1ac, \ b_2=bdb_1, \ w=a_1w_1b_1$ and
\begin{eqnarray*}
a_1s_1b_1-a_2s_2b_2&=&a_1s_1b_1-a_1acs_2bdb_1\\
&=&a_1(s_1-acs_2bd)b_1\\
&=&a_1(s_1,s_2)_{w_1}b_1\\
&\equiv&0 \ \ \ \ \ mod(S,w).
\end{eqnarray*}

\emph{(1.4)} $X,Y$-skew-inclusion

Assume that $w_1^X=\bar{s_1}^X=a\bar{s_2}^Xb, \ a,b\in X^*$ and
$w_1^Y=\bar{s_2}^Y=c\bar{s_1}^Yd, \ c,d\in Y^*$. Then
$a_2^X=a_1^Xa$, $b_2^X=bb_1^X$, $a_1^Y=a_2^Yc$ and $b_1^Y=db_2^Y$.
Thus, $w=a_1^Xa_2^Yw_1b_1^Xb_2^Y$ and
\begin{eqnarray*}
a_1s_1b_1-a_2s_2b_2&=&a_1^Xa_2^Ycs_1b_1^Xdb_2^Y-a_1^Xaa_2^Ys_2bb_1^Xb_2^Y\\
&=&a_1^Xa_2^Y(cs_1d-as_2b)b_1^Xb_2^Y\\
&=&a_1^Xa_2^Y(s_1,s_2)_{w_1}b_1^Xb_2^Y\\
&\equiv&0 \ \ \ \ \ mod(S,w). \end{eqnarray*}

\emph{Case 2} Intersection

\emph{(2.1)} $X$-intersection only

We may assume that $\bar{s_1}^X$ is at the left of $\bar{s_2}^X$,
i.e., $w_1^X=\bar{s_1}^Xb=a\bar{s_2}^X$, $a,b\in X^*$ and
$|\bar{s_1}^X|+|\bar{s_2}^X|>|w_1^X|$. Then $a_2^X=a_1^Xa$ and
$b_1^X=bb_2^X$. There are two cases to be consider:
$w_1^Y=\bar{s_1}^Yc\bar{s_2}^Y$ and $w_1^Y=\bar{s_2}^Yc\bar{s_1}^Y,
\ c\in Y^*.$

For $w_1^Y=\bar{s_1}^Yc\bar{s_2}^Y$, i.e.,
$w_1=\bar{s_1}bc\bar{s_2}^Y$, we have $a_2^Y=a_1^Y\bar{s_1}^Yc$,
$b_1^Y=c\bar{s_2}^Yb_2^Y$, $w=a_1\bar{s_1}ac\bar{s_2}b_2=a_1w_1b_2$
and
\begin{eqnarray*}
a_1s_1b_1-a_2s_2b_2&=&a_1s_1bb_2^Xc\bar{s_2}^Yb_2^Y-a_1^Xaa_1^Y\bar{s_2}^Ycs_2b_2\\
&=&a_1(s_1bc\bar{s_2}^Y-a\bar{s_2}^Ycs_2)b_2\\
&=&a_1(s_1,s_2)_{w_1}b_2\\
&\equiv&0 \ \ \ \ \ \ \ mod(S,w).
\end{eqnarray*}

For $w_1^Y=\bar{s_2}^Yc\bar{s_1}^Y$, i.e.,
$w_1=\bar{s_2}^Yc\bar{s_1}b$, we have $a_1^Y=a_2^Y\bar{s_2}^Yc$,
$b_2^Y=c\bar{s_1}^Yb_1^Y, \ w=a_1^Xa_2^Yw_1b_2^Xb_1^Y$ and
\begin{eqnarray*}
a_1s_1b_1-a_2s_2b_2&=&a_1^Xa_2^Y\bar{s_2}^Ycs_1bb_2^Xb_1^Y-a_1^Xaa_2^Ys_2b_2^Xc\bar{s_1}^Yb_1^Y\\
&=&a_1^Xa_2^Y(\bar{s_2}^Ycs_1b-as_2c\bar{s_1}^Y)b_2^Xb_1^Y\\
&=&a_1^Xa_2^Y(s_1,s_2)_{w_1}b_2^Xb_1^Y\\
&\equiv&0 \ \ \ \ \ \ \ mod(S,w).
\end{eqnarray*}

\emph{(2.2)} $Y$-intersection only

This case is similar to (2.1).

\emph{(2.3)} $X,Y$-intersection

Assume that $w_1^X=\bar{s_1}^Xb=a\bar{s_2}^X$,
$w_1^Y=\bar{s_1}^Yd=c\bar{s_2}^Y$, $a,b\in X^*, \ c,d\in Y^*$,
$|\bar{s_1}^X|+|\bar{s_2}^X|>|w_1^X|$ and
$|\bar{s_1}^Y|+|\bar{s_2}^Y|>|w_1^Y|$. Then $a_2^X=a_1^Xa$,
$b_1^X=bb_2^X$, $a_2^Y=a_1^Yc$,  $b_1^Y=db_2^Y$, $w=a_1w_1b_2$ and

\begin{eqnarray*}
a_1s_1b_1-a_2s_2b_2&=&a_1s_1bb_2^Xdb_2^Y-a_1^Xaa_1^Ycs_2b_2\\
&=&a_1(s_1bd-acs_2)b_2\\
&=&a_1(s_1,s_2)_{w_1}b_2\\
&\equiv&0 \ \ \ \ \ mod(S,w).
\end{eqnarray*}

\emph{(2.4)} $X,Y$-skew-intersection

Assume that $w_1^X=\bar{s_1}^Xb=a\bar{s_2}^X$,
$w_1^Y=c\bar{s_1}^Y=\bar{s_2}^Yd$,
$|\bar{s_1}^X|+|\bar{s_2}^X|>|w_1^X|$,
$|\bar{s_1}^Y|+|\bar{s_2}^Y|>|w_1^Y|$, $a,b\in X^*, \ c,d\in Y^*$.
Then $a_2^X=a_1^Xa$, $b_1^X=bb_2^X$, $a_1^Y=a_2^Yc$, $b_2^Y=db_1^Y$,
$w=a_1^Xa_2^Yw_1b_2^Xb_1^Y$ and
\begin{eqnarray*}
a_1s_1b_1-a_2s_2b_2&=&a_1^Xa_2^Ycs_1bb_2^Xb_1^Y-a_1^Xaa_2^Ys_2b_2^Xdb_1^Y\\
&=&a_1^Xa_2^Y(cs_1b-as_2d)b_2^Xb_1^Y\\
&=&a_1^Xa_2^Y(s_1,s_2)_{w_1}b_2^Xb_1^Y \\
&\equiv&0 \ \ \ \ \ \ \ mod(S,w).
\end{eqnarray*}

\emph{Case 3} Both inclusion and intersection

\emph{(3.1)} $X$-inclusion and $Y$-intersection

We may assume that $w_1^X=\bar{s_1}^X=a\bar{s_2}^Xb$,  $a,b\in X^*$.
Then $a_2^X=a_1^Xa$ and $b_2^X=bb_1^X$. There two cases to consider:
$w_1^Y=\bar{s_1}^Yd=c\bar{s_2}^Y$ and
$w_1^Y=c\bar{s_1}^Y=\bar{s_2}^Yd$, where $c,d\in Y^*$,
$|\bar{s_1}^Y|+|\bar{s_2}^Y|>|w_1^Y|$.

For  $w_1^Y=\bar{s_1}^Yd=c\bar{s_2}^Y$, we have $a_2^Y=a_1^Yc$,
$b_1^Y=db_2^Y, \ w=a_1w_1b_1^Xb_2^Y$ and
\begin{eqnarray*}
a_1s_1b_1-a_2s_2b_2&=&a_1s_1b_1^Xdb_2^Y-a_1^Xaa_1^Ycs_2bb_1^Xb_2^Y\\
&=&a_1(s_1d-acs_2b)b_1^Xb_2^Y\\
&=&a_1(s_1,s_2)_{w_1}b_2\\
&\equiv&0 \ \ \ \ \ mod(S,w).
\end{eqnarray*}

For $w_1^Y=c\bar{s_1}^Y=\bar{s_2}^Yd$, we have $a_1^Y=a_2^Yc$,
$b_2^Y=db_1^Y$, $w=a_1^Xa_2^Yw_1b_2^Xdb_1^Y$ and
\begin{eqnarray*}
a_1s_1b_1-a_2s_2b_2&=&a_1^Xa_2^Ycs_1b_1-a_1^Xaa_2^Ys_2bb_1^Xdb_1^Y\\
&=&a_1^Xa_2^Y(cs_1-as_2bd)b_1\\
&=&a_1^Xa_2^Y(s_1,s_2)_{w_1}b_1 \\
&\equiv&0 \ \ \ \ \ \ \ mod(S,w).
\end{eqnarray*}

\emph{(3.2)} $X$-intersection and $Y$-inclusion

 Assume that
$w_1^X=\bar{s_1}^Xb=a\bar{s_2}^X$, $a,b\in Y^*$ with
$|\bar{s_1}^X|+|\bar{s_2}^X|>|w_1^X|$. Then $a_2^X=a_1^Xa$,
$b_1^X=bb_2^X$. There are  two cases to consider:
$w_1^Y=\bar{s_1}^Y=c\bar{s_2}^Yd$ and $\bar{s_2}^Y=c\bar{s_1}^Yd$,
where $c,d\in Y^*$.

For  $w_1^Y=\bar{s_1}^Y=c\bar{s_2}^Yd$, we have $a_2^Y=a_1^Yc$,
$b_2^Y=db_1^Y$, $w=a_1w_1b_2^Xb_1^Y$ and
\begin{eqnarray*}
a_1s_1b_1-a_2s_2b_2&=&a_1s_1bb_2^Xb_1^Y-a_1acs_2b_2^Xdb_1^Y\\
&=&a_1(s_1b-acs_2d)b_2^Xb_1^Y\\
&=&a_1(s_1,s_2)_{w_1}b_2^Xb_1^Y\\
&\equiv&0 \ \ \ \ \ mod(S,w).
\end{eqnarray*}

For $w_1^Y=\bar{s_2}^Y=c\bar{s_1}^Yd$, we have $a_1^Y=a_2^Yc$,
$b_1^Y=db_2^Y$, $w=a_1^Xa_2^Yw_1b_2$ and
\begin{eqnarray*}
a_1s_1b_1-a_2s_2b_2&=&a_1^Xa_2^Ycs_1bb_2^Xdb_2^Y-a_1^Xaa_2^Ys_2b_2\\
&=&a_1^Xa_2^Y(cs_1bd-as_2)b_2\\
&=&a_1^Xa_2^Y(s_1,s_2)_{w_1}b_2\\
&\equiv&0 \ \ \ \ \ mod(S,w).
\end{eqnarray*}

\emph{Case 4}. $\bar{s_1}$ and $\bar{s_2}$
  disjoint

For $w=w^Xw^Y$, by symmetry, there are two cases to consider:
$w^Y=a_1^Y\bar{s_1}^Yc\bar{s_2}^Yb_2^Y$ and
$w^Y=a_2^Y\bar{s_2}^Yc\bar{s_1}^Yb_1^Y$, where
$w^X=a_1^X\bar{s_1}^Xa\bar{s_2}^Xb_2^X, \ a\in X^*,  \
a_2^X=a_1^X\bar{s_1}^Xa, \ b_1^X=a\bar{s_2}^Xb_2^X$ and $c\in Y^*$.

For
$w=a_1^X\bar{s_1}^Xa\bar{s_2}^Xb_2^Xa_1^Y\bar{s_1}^Yc\bar{s_2}^Yb_2^Y=a_1\bar{s_1}ac\bar{s_2}b_2$,
 we have $a_2=a_1\bar{s_1}ac, \ b_1=ac\bar{s_2}b_2$ and
\begin{eqnarray*}
a_1s_1b_1-a_2s_2b_2&=&a_1s_1ac\bar{s_2}b_2-a_1\bar{s_1}acs_2b_2\\
&=&a_1(s_1-\bar{s_1})acs_2b_2-a_1s_1ac(s_2-\bar{s_2})b_2\\
&\equiv&0 \ \ \ \ \ \ \ mod(S,w).
\end{eqnarray*}

For
$w=a_1^X\bar{s_1}^Xa\bar{s_2}^Xb_2^Xa_2^Y\bar{s_2}^Yc\bar{s_1}^Yb_1^Y$,
we have $a_1^Y=a_2^Y\bar{s_2}^Yc, \ b_2^Y=c\bar{s_1}^Yb_1^Y$ and

\begin{eqnarray*}
a_1s_1b_1-a_2s_2b_2&=&a_1^Xa_2^Y\bar{s_2}^Ycs_1a\bar{s_2}^Xb_2^Xb_1^Y-a_1^X\bar{s_1}^Xaa_2^Ys_2b_2^Xc\bar{s_1}^Yb_1^Y\\
&=&a_1^Xa_2^Y(\bar{s_2}^Ycs_1a\bar{s_2}^X-\bar{s_1}^Xas_2c\bar{s_1}^Y)b_2^Xb_1^Y.\\
\end{eqnarray*}
Let $s_1=\sum_{i=1}^n\alpha_iu_{1i}^Xu_{1i}^Y$ and
$s_2=\sum_{j=1}^m\beta_ju_{2j}^Xu_{2j}^Y$, where
$\alpha_1=\beta_1=1.$ Then
\begin{eqnarray*}
\bar{s_2}^Ycs_1a\bar{s_2}^X-\bar{s_1}^Xas_2c\bar{s_1}^Y&=&
\sum_{i=2}^n\alpha_iu_{1i}^Xa\bar{s_2}cu_{1i}^Y-\sum_{j=2}^m\beta_iu_{2j}^Yc\bar{s_1}au_{2j}^X\\
&=&\sum_{i=2}^n\alpha_iu_{1i}^Xa(\bar{s_2}-s_2)cu_{1i}^Y+\sum_{j=2}^m\beta_ju_{2j}^Yc(s_1-\bar{s_1})au_{2j}^X\\
& \ &+\sum_{i=2}^n\alpha_iu_{1i}^Xas_2cu_{1i}^Y-\sum_{j=2}^m\beta_ju_{2j}^Ycs_1au_{2j}^X\\
&\equiv&\sum_{i=2}^n\sum_{j=2}^m\alpha_i\beta_ju_{1i}^Xau_{2j}^Xu_{2j}^Ycu_{1i}^Y-
\sum_{j=2}^m\sum_{i=2}^n\alpha_i\beta_ju_{2j}^Ycu_{1i}^Yu_{1i}^Xau_{2j}^X\\
&\equiv&0 \ \  \ \ \ \ \ mod(S,w_1),
\end{eqnarray*}
where
$w_1=\bar{s_2}^Yc\bar{s_1}a\bar{s_2}^X=\bar{s_1}^Xa\bar{s_2}c\bar{s_1}^Y$.
Since $w=a_1^Xa_2^Yw_1b_2^Xb_1^Y$, we have
\begin{eqnarray*}
a_1s_1b_1-a_2s_2b_2&=&a_1^Xa_2^Y(\bar{s_2}^Ycs_1a\bar{s_2}^X-\bar{s_1}^Xas_2c\bar{s_1}^Y)b_2^Xb_1^Y\\
&\equiv&0 \ \  \ \ \ \ \ mod(S,w).
\end{eqnarray*}
This completes the proof.\ \ $\square$

\begin{lemma}\label{2.6}
Let $S\subset k\langle X\rangle \otimes k\langle Y\rangle$ with each
$s\in S$ monic and $Irr(S)=\{w \in N|w\neq a\overline{s}b, \ a,b\in
N, \ s\in S\}$. Then for any $f\in k\langle X\rangle \otimes
k\langle Y\rangle$,
$$
f=\sum_{a_i \bar{s_i}b_i\leq \bar{f}} \alpha_i a_i
s_ib_i+\sum_{u_j\leq \bar{f}} \beta_j u_j,
$$
where $\alpha_i, \ \beta_j \in k, \ a_i,b_i\in N, \ s_i\in S \
\mbox{and} \ u_j\in Irr(S)$.
\end{lemma}
{\bf Proof.} Let  $f=\sum\limits_{i}\alpha_{i}u_{i}\in k\langle
X\rangle \otimes k\langle Y\rangle$, where $0\neq{\alpha_{i}\in{k}}$
and $u_{1}>u_{2}>\cdots$. If $u_1\in{Irr(S)}$, then let
$f_{1}=f-\alpha_{1}u_1$. If $u_1\not\in{Irr(S)}$, then there exist
some $s\in{S}$ and $a_1,b_1\in N$, such that $\bar
f=a_1\bar{s_1}b_1$. Let $f_1=f-\alpha_1a_1s_1b_1$. In both cases, we
have $\bar{f_1}<\bar{f}$. Then the result follows from the induction
on $\bar{f}$. \ \ \ \ $\square$

From the above lemmas, we reach the following theorem:

\begin{theorem}
(Composition-Diamond lemma for tensor product $k\langle
X\rangle\otimes k\langle Y\rangle$)\label{t1} Let $S\subset k\langle
X\rangle\otimes k\langle Y\rangle$ with each $s\in S$ monic and $<$
the order on $N=X^*Y^*$ as before. Then the following statements are
equivalent:
\begin{enumerate}
\item[(1)]\ $S$ is a Gr\"obner-Shirshov basis in $k\langle X\cup Y|T\rangle$.
\item[(2)]\
$f\in Id(S)\Rightarrow\overline{f}=a\overline{s}b$ for some $a,b\in
N , \ s\in S$.
\item[(3)]\
$Irr(S)=\{w \in N|w\neq a\overline{s}b, \ a,b\in N, \ s\in S\}$ is a
$k$-linear basis for the factor $k\langle X\cup Y|T\rangle/Id(S)$.
\end{enumerate}
\end{theorem}

\textbf{Proof: }$(1)\Rightarrow (2)$. Suppose that $0\neq f\in
Id(S)$. Then $f=\sum \alpha_i a_i s_i b_i$ for some $\alpha_i\in k,
\ a_i,b_i\in N, \ s_i\in S$. Let $w_i=a_i \overline{s}_ib_i$ and
$w_1=w_2=\cdots=w_l>w_{l+1}\geq \cdots$. We will prove that
$\overline{f}=a\overline{s}b$ for some $a,b\in N, \ s\in S$, by
using induction on $l$ and $w_1$. If $l=1$, then the result is
clear. If $l>1$, then $w_1=a_1 \bar{s_1}b_1=a_2 \bar{s_2}b_2$. Now,
by (1) and Lemma \ref{l3}, $a_1s_1b_1\equiv a_2s_2b_2 \ \
mod(S,w_1)$. Thus,
\begin{eqnarray*}
\alpha_1a_1 s_1b_1+\alpha_2a_2 s_2
b_2&=&(\alpha_1+\alpha_2)a_1s_1b_1+\alpha_2(a_2
s_2b_2-a_1s_1b_1)\\
&\equiv&(\alpha_1+\alpha_2)a_1s_1b_1 \ \ \ \ \ \ mod(S,w_1).
\end{eqnarray*}
 By induction on $l$ and $w_1$, we have the
result.

$(2)\Rightarrow (3)$. For any $0\neq f\in k\langle X\cup
Y|T\rangle$, by Lemma \ref{2.6}, we can express $f$ as
$$
f=\sum \alpha_i a_i s_ib_i+\sum \beta_j u_j,
$$
where $\alpha_i, \ \beta_j \in k, \ a_i,b_i\in N, \ s_i\in S \
\mbox{and} \ u_j\in Irr(S)$. Then $Irr(S)$ generates the factor
algebra. Moreover, if $0\neq h=\sum \beta_j u_j\in Id(S)$, $u_j\in
Irr(S), u_1>u_2>\cdots  \ and  \ \beta_1\neq 0$, then
$u_1=\bar{h}=a\bar{s}b$ for some $a,b\in N, \ s\in S$ by (2), a
contradiction. This shows that $Irr(S)$ is a linear basis of the
factor algebra.

$(3)\Rightarrow (1)$. For any $f, \ g\in S$, we have $h=(f,g)_w\in
Id(S)$. The result is trivial if $(f,g)_w=0$. Assume that
$(f,g)_w\neq 0$. Then, by Lemma \ref{2.6} and (3), we have
$$
h=\sum_{a_i \bar{s_i}b_i\leq \bar{h}} \alpha_i a_i s_ib_i.
$$
Now, by noting that $\bar{h}=\overline{(f,g)_w}<w$, we know that (1)
holds. \ \ $\square$

\ \

\noindent{\bf Remark}: Theorem \ref{t1} is valid for any monomial
order on $X^*Y^*$.

\ \

\noindent{\bf Remark}: Theorem \ref{t1} is exact the
Composition-Diamond lemma for associative algebras (Lemma
\ref{l1.4}) when $Y=\emptyset$.

\section{Applications}

Now, we give some applications of Theorem \ref{t1}.

\begin{example}\label{ex1} Suppose that for the deg-lex order,
$S_1$ and $S_2$ are  Gr\"{o}bner-Shirshov bases in $k\langle
X\rangle$ and $k\langle Y\rangle$ respectively. Then  for the
deg-lex order on $X^*Y^*$, $S_1\cup S_2$ is a Gr\"{o}bner-Shirshov
basis in $k\langle X\cup Y|T\rangle=k\langle X\rangle\otimes
k\langle Y\rangle$. It follows that  $k\langle X|S_1\rangle\otimes
k\langle Y|S_2\rangle=k\langle X\cup Y|T\cup S_1\cup S_2 \rangle$.
\end{example}

\textbf{Proof: } The possible compositions in $S_1\cup S_2$ are
$X$-including only, $X$-intersection only, $Y$-including only and
$Y$-intersection only. Suppose $f,g\in S_1$ and $(f,g)_{w_1}\equiv 0
\ \ mod(S_1,w_1)$ in $k\langle X\rangle$. Then in $k\langle X\cup
Y|T\rangle$, $(f,g)_w=(f,g)_{w_1}c$, where $w=w_1c$ for any $c\in
Y^*$. From this it follows that each composition in $S_1\cup S_2$ is
trivial modulo $S_1\cup S_2$. \ \ $\square$

A special case of Example \ref{ex1} is the following.
\begin{example} Let $X,Y$ be linearly ordered sets, $k[X]$ the free
commutative  associative algebra generated by $X$. Then
$S=\{x_ix_j=x_jx_i|x_i>x_j, x_i,x_j\in X \}$ is a
Gr\"{o}bner-Shirshov basis in $k\langle X\rangle\otimes k\langle
Y\rangle$ with respect to the deg-lex order. Therefore, $k[X]\otimes
k\langle Y \rangle=k\langle X\cup Y|T\cup S\rangle$.
\end{example}

In \cite{eisenbud}, a Gr\"{o}bner-Shirshov basis in $k\langle
X\rangle$ is constructed by lifting a commutative Gr\"{o}bner basis
and adding commutators. Let $X$ be a well-ordered set, $[X]$ the
free commutative monoid generated by $X$ and $k[X]$ the polynomial
ring. Let $S_1=\{h_{ij}=x_ix_j-x_jx_i| \ i>j\}\subset k\langle X
\rangle$. Consider the natural map $\gamma: k\langle X
\rangle\rightarrow k[X] $ taking $x_i$ to $x_i$ and the
\emph{lexicographic splitting } of $\gamma$, which is defined as the
$k$-linear map

$$\delta : k[X]\rightarrow k\langle X \rangle, \ \
x_{i_1}x_{i_2}\cdots x_{i_r}\mapsto x_{i_1}x_{i_2}\cdots x_{i_r} \ \
if \ \ i_1\leq i_2\cdots \leq i_r.$$

For any $u\in [X]$, we present $u=x_1^{l_1}x_2^{l_2}\cdots
x_n^{l_n}$, where $l_i\geq 0$. We  use any monomial order on $[X]$.
For any $f\in k[X]$, $\bar{f}$ means the leading monomial of $f$.

Following \cite{eisenbud}, we define an order on $X^*$ using the
order $x_1<x_2<\cdots<x_n$ as follows: for any $u,v\in X^*$,
$$
u>v\Leftrightarrow \gamma(u)>\gamma(v) \ \mbox{ in } \ [X] \ \mbox{
or } \ (\gamma(u)=\gamma(v) \ and \ u>_{lex}v).
$$
It is easy to check that this order is monomial  on $X^*$ and
$\overline{\delta(s)}=\delta(\bar{s})$ where $s\in k[X]$.
Moreover, for any $v\in \gamma^{-1}(u)$, $v\geq\delta(u)$.

 For any $m=x_{i_1}x_{i_2}\cdots x_{i_r}\in [X], \ i_1\leq i_2\cdots
\leq i_r$, denote the set of all the monomials $u\in
[x_{i_1+1},\cdots, x_{i_r-1}]$ by $U(m)$.

The proofs of the following lemmas are straightforward.

\begin{lemma}\label{4.1} Let $a,b\in X^*, \ a=\delta(\gamma(a)), \
b=\delta(\gamma(b))$ and $ s\in k[X]$. If
$w=a\delta(\bar{s})b=\delta(\gamma(ab)\bar{s})$, then, in $k\langle
X\rangle$,
$$a\delta(s)b\equiv\delta(\gamma(ab)s) \ \ \ \ mod(S_1,w).$$
\end{lemma}

\textbf{Proof: }Suppose that $s=\bar{s}+s'$ and
$h=a\delta(s)b-\delta(\gamma(ab)s)$. Since
$a\delta(\bar{s})b=\delta(\gamma(ab)\bar{s})$, we  have
$h=a\delta(s')b-\delta(\gamma(ab)s')$, and $\bar{h}<w$. By noting
that $\gamma(\delta(\gamma(ab)s')=\gamma(\delta(\gamma(ab)s'))$,
$h\equiv0 \ \ \ mod(S_1,w)$. \ \ $\square$

\begin{lemma}\label{4.2} Let $f,g\in k[X], \ \bar{g}=x_{i_1}x_{i_2}\cdots
x_{i_r} \ (i_1\leq i_2\leq\cdots\leq i_r ) \ and \
w=\delta(\bar{f}\bar{g}).$ Then, in $k\langle X\rangle$,
$$\delta((f-\bar{f})g)\equiv
\sum\alpha_ia_i\delta(u_ig)b_i \ \ \ mod(S_1,w)$$
 where $\alpha_i\in
k, \ a_i\in [x\in X  |  x\leq x_{i_1}], \ b_i\in [x\in X  |  x\geq
x_{i_r}], \ u_i\in U(\bar{g})$ and
$\gamma(\sum\alpha_ia_iu_ib_i)=f-\bar{f}$.
\end{lemma}

\begin{theorem}\label{th9}(\cite{eisenbud}) Let the orders on $[X]$ and $X^*$
be defined as above. If $S$ is a minimal Gr\"{o}bner basis in
$k[X]$, then $S'=\{\delta(us)  |  s\in S, u\in U(\bar{s}) \} \cup
S_1$ is a Gr\"{o}bner-Shirshov basis in $k\langle X\rangle$.
\end{theorem}

\textbf{Proof: }We will show that all the possible compositions of
elements in $S'$ are trivial. Let $f=\delta(us_1), \ g=\delta(vs_2)$
and $h_{ij}=x_ix_j-x_jx_i \in S'$.

$(1)$  $f \wedge g$

Case 1. $f$ and $g$ have a composition of including, i.e.,
$w=\delta(u\bar{s_1})=a\delta(v\bar{s_2})b$ for some $a,b\in X^*$
and $a=\delta(\gamma(a)),b=\delta(\gamma(b))$.

If  $s_1$ and $s_2$ have no composition in $k[X]$, i.e.,
$lcm(\bar{s_1}\bar{s_2})=\bar{s_1}\bar{s_2}$, then $u=u'\bar{s_2}, \
\gamma(ab)v=u'\bar{s_1}$ for some $u'\in [X]$. By Lemma \ref{4.1}
and Lemma \ref{4.2}, we have
\begin{eqnarray*}
(f,g)_w &=& \delta(us_1)-a\delta(vs_2)b \\
&\equiv&\delta(us_1)-\delta(\gamma(ab)vs_2) \\
&\equiv&\delta(u'\bar{s_2}s_1)-\delta(u'\bar{s_1}s_2) \\
&\equiv&\delta(u'(s_1-\bar{s_1})s_2)-\delta(u'(s_2-\bar{s_2})s_1)  \\
&\equiv&0 \ \ \ \ \ \ \ mod(S',w).
 \end{eqnarray*}

Since, in $k[X]$, $S$ is a minimal Gr\"{o}bner basis, the possible
compositions are only intersection. If $s_1$ and $s_2$ have
composition of intersection in $k[X]$, i.e.,
$(s_1,s_2)_{w'}=a's_1-b's_2$, where $a',b'\in[X], \
w'=a'\bar{s_1}=b'\bar{s_2}$ and $|w'|<|\bar{s_1}|+|\bar{s_2}|$, then
$w'$ is a subword of $\gamma(w)$.
 Therefore, we have
 $w=\delta(tw')=\delta(ta'\bar{s_1})=\delta(tb'\bar{s_2})$ and $u=ta',
 \gamma(ab)v=tb'$ for some $t\in [X]$. Then
\begin{eqnarray*}
(f,g)_w &=& \delta(us_1)-a\delta(vs_2)b \\
&\equiv&\delta(us_1)-\delta(\gamma(ab)vs_2) \\
&\equiv&\delta(ta's_1)-\delta(tb's_2) \\
&\equiv&\delta(t(a's_1-b's_2))  \\
&\equiv&\delta(t(s_1,s_2)_{w'})\\
&\equiv&0 \ \ \ \ \ \ \ mod(S',w)
 \end{eqnarray*}
since $t\overline{(s_1,s_2)_{w'}}<tw'=\gamma(w).$

Case 2. If $f$ and $g$ have a composition of intersection, we may
assume that $\bar{f}$ is on the left of $\bar{g}$,  i.e.,
$w=\delta(u\bar{s_1})a=b\delta(v\bar{s_2})$ for some $a,b\in X^*$
and $a=\delta\gamma(a), b=\delta\gamma(b)$. Similarly to Case 1, we
have to consider whether $s_1$ and $s_2$ have compositions in $k[X]$
or not. One can check that both cases are trivial mod$(S',w)$ by
Lemma \ref{4.1} and Lemma \ref{4.2}.

$(2)$  $f \wedge h_{ij}$

By noting that $\overline{h_{ij}}=x_ix_j$ can not be a subword of
$\bar{f}=\delta(u\bar{s_1})$ since $i>j$,  only possible
compositions are intersection. Suppose that $\bar{s_1}=x_{i_1}\cdots
x_{i_r}x_i, \ (i_1\leq i_2\leq\cdots\leq i_r\leq i)$. Then
$\bar{f}=\delta(u\bar{s_1})=x_{i_1}vx_i$ for some $v\in k\langle
X\rangle,  \ v=\delta\gamma(v)$ and $w=\delta(u\bar{s_1})x_j$.

If $j\leq i_1$, then
\begin{eqnarray*}
(f,h_{ij})_w &=& \delta(us_1)x_j- x_{i_1}v(x_ix_j-x_jx_i)\\
&=&\delta(u(s_1-\bar{s_1}))x_j+ x_{i_1}vx_jx_i \\
&\equiv& x_j\delta(u(s_1-\bar{s_1}))+x_j x_{i_1}vx_i \\
&\equiv& x_j(\delta(u(s_1-\bar{s_1}))+\delta(u\bar{s_1}) ) \\
&\equiv&x_j\delta(us_1)\\
&\equiv&0 \ \ \ \ \ \ \ mod(S',w).
 \end{eqnarray*}

If $j> i_1$, then $ux_j\in U(\bar{s_1})$ and
\begin{eqnarray*}
(f,h_{ij})_w &=& \delta(us_1)x_j- x_{i_1}v(x_ix_j-x_jx_i)\\
&=&\delta(u(s_1-\bar{s_1}))x_j+ x_{i_1}vx_jx_i \\
&\equiv& \delta(ux_j(s_1-\bar{s_1}))+\delta( x_{i_1}vx_ix_j) \\
&\equiv& \delta(ux_j(s_1-\bar{s_1}))+\delta(ux_j\bar{s_1}) \\
&\equiv&\delta(ux_js_1)\\
&\equiv&0 \ \ \ \ \ \ \ mod(S',w).
 \end{eqnarray*}

 Then we complete the proof. \ \ $\square$

Now we extend $\gamma$  and  $\delta$  as follows.
\begin{eqnarray*}
\gamma\otimes\mathbf{1}: \  k\langle X\rangle\otimes k\langle
Y\rangle&\rightarrow& k[X]\otimes k\langle Y\rangle, \
u^Xu^Y\mapsto\gamma(u^X)u^Y,
\\
\delta\otimes\mathbf{1}: \ k[X]\otimes k\langle Y\rangle
&\rightarrow& k\langle X\rangle\otimes k\langle Y\rangle, \
u^Xu^Y\mapsto \delta(u^X)u^Y.
\end{eqnarray*}

Any polynomial $f\in k[X]\otimes k\langle Y \rangle$ has a
presentation $f=\sum \alpha_iu^X_iu^Y_i$, where $\alpha_i\in
k,u^X_i\in [X] \ and \ u^Y_i\in Y^*$.

Let the orders on $[X]$ and $Y^*$ be any monomial oeders
respectively. We order the set $[X]Y^*=\{u=u^Xu^Y|u^X\in [X], \
u^Y\in Y^*\}$ as follows. For any $u,v\in [X]Y^*$,
$$
u>v\Leftrightarrow u^Y>v^Y \ or \ (u^Y=v^Y  \ \mbox{and} \  \
u^X>v^X).
$$

 Now, we order $X^*Y^*$: for any $u,v\in X^*Y^*$,
 $$
 u>v\Leftrightarrow
\gamma(u^X)u^Y>\gamma(v^X)v^Y \ \mbox{or} \
(\gamma(u^X)u^Y=\gamma(v^X)v^Y \ \mbox{and} \ u^X>_{lex}v^X) .
$$
This order is clearly a monomial order on $X^*Y^*$.

The following definitions of compositions and Gr\"{o}bner-Shirshov
bases are essentially from \cite{MZ}.

Let $f,g$ be monic polynomials of $k[X]\otimes k\langle Y \rangle$,
 $L$ the least common multiple of $\bar{f}^X$ and $\bar{g}^X$.

$1.$  Inclusion

Let $\bar{g}^Y$ be a subword of $\bar{f}^Y$, say,
$\bar{f}^Y=c\bar{g}^Yd$ for some $c,d\in Y^*$. If
$\bar{f}^Y=\bar{g}^Y$ then $\bar{f}^X\geq \bar{g}^X$ and if
$\bar{g}^Y=1$ then we set $c=1$. Let $w=L\bar{f}^Y=Lc\bar{g}^Yd$. We
define the composition

$$C_1(f,g,c)_w= \frac{L}{\bar{f}^X}f-\frac{L}{\bar{g}^X}cgd.$$

$2.$  Overlap

Let a non-empty beginning of $\bar{g}^Y$ be a non-empty ending of
$\bar{f}^Y$, say,
$\bar{f}^Y=cc_0,\bar{g}^Y=c_0d,\bar{f}^Yd=c\bar{g}^Y$ for some
$c,d,c_0\in Y^*$ and $c_0\neq1$. Let $w=L\bar{f}^Yd=Lc\bar{g}^Y$. We
define the composition
$$C_2(f,g,c_0)_w= \frac{L}{\bar{f}^X}fd-\frac{L}{\bar{g}^X}cg.$$

$3.$  External

Let $c_0\in Y^*$ be any associative word (possibly empty). In the
case that the greatest common divisor of $\bar{f}^X$ and $\bar{g}^X$
is non-empty and $\bar{f}^Y, \bar{g}^Y$ are non-empty, we define the
composition
$$
C_3(f,g,c_0)_w=\frac{L}{\bar{f}^X}fc_0\bar{g}^Y-\frac{L}{\bar{g}^X}\bar{f}^Yc_0g,
$$
where $w=L\bar{f}^Yc_0\bar{g}^Y$.

Let $S$ be a monic subset of $k[X]\otimes k\langle Y \rangle$. Then
$S$  is called a Gr\"{o}bner-Shirshov basis (standard basis) if for
any element $f\in Id(S)$, $\bar{f}$ contains $\bar{s}$ as its
subword for some $s\in S$.

It is defined as usual that a composition is trivial modulo $S$ and
corresponding $w$. We also have that $S$ is a Gr\"{o}bner-Shirshov
basis in $k[X]\otimes k\langle Y \rangle$ if and only if all the
possible compositions of its elements are trivial. A
Gr\"{o}bner-Shirshov basis in $k[X]\otimes k\langle Y \rangle$ is
called minimal if for any $s\in S$ and all $s_i\in S\setminus
\{s\}$, $\bar{s_i}$ is not a subword of $\bar{s}$.

Similar to the proof of Theorem \ref{th9}, we have the following
theorem.

\begin{theorem}  Let the orders on $[X]Y^*$ and $X^*Y^*$ be defined as before.
If $S$ is a minimal Gr\"{o}bner-Shirshov  basis in $k[X]\otimes
k\langle Y \rangle$, then  $S'=\{\delta(us) | s\in S, u\in
U(\bar{s}^X) \} \cup S_1$ is a Gr\"{o}bner-Shirshov basis in
$k\langle X\rangle\otimes k\langle Y\rangle$, where
$S_1=\{h_{ij}=x_ix_j-x_jx_i| \ i>j\}$.
\end{theorem}

\textbf{Proof: }We will show that all the possible compositions of
elements in $S'$ are trivial.

For $s_1,s_2\in S$, let $f=\delta(us_1), \ g=\delta(vs_2), \
h_{ij}=x_ix_j-x_jx_i \in S'$ and $L=lcm(\bar {s_1}^X,\bar {s_2}^X)$.

$1.$  $f \wedge g$

In this case, all the possible compositions of $f \wedge g$  are
related the ambiguities $w$'s (in the following,  $a,b\in X^*, \
c,d\in Y^*$).

\emph{(1.1)} $X$-inclusion only

$w^X=\delta(u\bar {s_1}^X)=a\delta(v\bar{s_2}^X)b$,
$w^Y=\bar{s_1}^Yc\bar{s_2}^Y$ or $w^Y=\bar{s_2}^Yc\bar{s_1}^Y$.

\emph{(1.2)} $Y$-inclusion only

$w^X=\delta(u\bar {s_1}^X)a\delta(v\bar{s_2}^X)$ or
$w^X=\delta(v\bar{s_2}^X)a\delta(u\bar {s_1}^X)$,
$w^Y=\bar{s_1}^Y=c\bar{s_2}^Yd$.

\emph{(1.3)} $X,Y$-inclusion

$w=\delta(u\bar{s_1})=ac\delta(v\bar{s_2})bd$.

\emph{(1.4)} $X,Y$-skew-inclusion

$w^X=\delta(u\bar{s_1})^X=a\delta(v\bar{s_2}^X)b$,
$w^Y=\bar{s_2}^Y=c\bar{s_1}^Yd$.

\emph{(2.1)} $X$-intersection only

$w^X=\delta(u\bar{s_1})^Xa=b\delta(v\bar{s_2}^X)$,
$w^Y=\bar{s_1}^Yc\bar{s_2}^Y$ or $w^Y=\bar{s_2}^Yc\bar{s_1}^Y$.

\emph{(2.2)} $Y$-intersection only

$w^X=\delta(u\bar{s_1})^Xa\delta(v\bar{s_2}^X)$ or
$w^X=\delta(u\bar{s_2})^Xa\delta(v\bar{s_1}^X)$,
$w^Y=\bar{s_1}^Yc=d\bar{s_2}^Y$.

\emph{(2.3)} $X,Y$-intersection

$w^X=\delta(u\bar{s_1})^Xa=b\delta(v\bar{s_2}^X)$,
$w^Y=\bar{s_1}^Yc=d\bar{s_2}^Y$.

\emph{(2.4)} $X,Y$-skew-intersection

$w^X=\delta(u\bar{s_1})^Xa=b\delta(v\bar{s_2}^X)$,
$w^Y=c\bar{s_1}^Y=\bar{s_2}^Yd$.

\emph{(3.1)} $X$-inclusion and $Y$-intersection

$w^X=\delta(u\bar{s_1})^X=a\delta(v\bar{s_2}^X)b$,
$w^Y=\bar{s_1}^Yc=d\bar{s_2}^Y$ or $w^Y=c\bar{s_1}^Y=\bar{s_2}^Yd$.

\emph{(3.2)} $X$-intersection and $Y$-inclusion

$w^X=\delta(u\bar{s_1})^Xa=b\delta(v\bar{s_2}^X)$,
$w^Y=\bar{s_1}^Y=c\bar{s_2}^Yd$ or $w^Y=\bar{s_2}^Y=c\bar{s_1}^Yd$.

We only check the  cases of (1.1), (1.2) and (1.3). Other cases are
similarly checked.

 \emph{(1.1)} $X$-inclusion only

Suppose that $w^X=\delta(u\bar {s_1}^X)=a\delta(v\bar{s_2}^X)b, \
a,b\in X^*$ and $\bar {s_1}^Y, \ \bar{s_2}^Y$ are disjoint.  There
are two cases to consider: $w^Y=\bar{s_1}^Yc\bar{s_2}^Y$ and
$w^Y=\bar{s_2}^Yc\bar{s_1}^Y$, where $c\in Y^*$. We  will only prove
the first case and the second is similar.

If  $s_1$ and $s_2$ have no composition in $k[X]\otimes k\langle Y
\rangle$, i.e., $lcm(\bar{s_1},\bar{s_2})=\bar{s_1}\bar{s_2}$, then
$u=u'\bar{s_2^X}, \ \gamma(ab)v=u'\bar{s_1^X}$ for some $u'\in [X]$.
By the proof of Theorem \ref{th9}, we have
\begin{eqnarray*}
(f,g)_w &=& \delta(us_1)c\bar{s_2}^Y-\bar{s_1}^Yca\delta(vs_2)b \\
&\equiv&\delta(us_1\gamma(c\bar{s_2}^Y))-\delta(\gamma(\bar{s_1}^Yc)\gamma(ab)vs_2) \\
&\equiv&\delta(u'\bar{s_2}^Xs_1c\bar{s_2}^Y)-\delta(\bar{s_1}^Ycu'\bar{s_1}^Xs_2) \\
&\equiv&\delta(u's_1c\bar{s_2})-\delta(u'\bar{s_1}cs_2) \\
&\equiv&\delta(u'(s_1-\bar{s_1})cs_2)-\delta(u's_1c(s_2-\bar{s_2}))  \\
&\equiv&0 \ \ \ \ \ \ \ mod(S',w).
 \end{eqnarray*}

If $s_1$ and $s_2$ have composition of external (the elements of $S$
have no composition of inclusion because $S$ is minimal and $s_1$
and $s_2$ have no composition of overlap because $s_1^Y$ and $s_2^Y$
 are disjoint ) in $k[X]\otimes k\langle Y \rangle$, i.e.,
$C_3(s_1,s_2,c)_{w'}=\frac{L}{\bar{s_1}^X}s_1\gamma(c\bar{s_2}^Y)-\frac{L}{\bar{s_2}^X}\gamma(\bar{s_1}^Yc)s_2
=t_2s_1\gamma(c\bar{s_2}^Y)-t_1\gamma(\bar{s_1}^Yc)s_2$ where
$gcd(\bar{s_1}^X,\bar{s_2}^X)=t\neq1,\bar{s_1}^X=tt_1,\bar{s_2}^X=tt_2$
and $L=tt_1t_2, \ w'=L\gamma(\bar{s_1}^Yc\bar{s_2}^Y)$, then $w'$ is
a subword of $\gamma(w)$.
 Therefore, we have
 $w=\delta(mw')$ and $u=mt_2,\gamma(ab)v=mt_1$ since
  $ut_1=\gamma(ab)vt_2$ and $gcd(t_1,t_2)=1$. Then
\begin{eqnarray*}
(f,g)_w &=& \delta(us_1)c\bar{s_2}^Y-\bar{s_1}^Yca\delta(vs_2)b \\
&\equiv&\delta(us_1\gamma(c\bar{s_2}^Y))-\delta(\gamma(\bar{s_1}^Yc)\gamma(ab)vs_2) \\
&\equiv&\delta(mt_2s_1\gamma(c\bar{s_2}^Y))-\delta(mt_1\gamma(\bar{s_1}^Yc)s_2) \\
&\equiv&\delta(mC_3(s_1,s_2,c)_{w'}) \\
&\equiv&0 \ \ \ \ \ \ \ mod(S',w)
 \end{eqnarray*}
since $m\overline{C_3(s_1,s_2,c)_{w'}}<mw'=\gamma(w).$

\emph{(1.2)} $Y$-inclusion only

Suppose that $w^Y=\bar{s_1}^Y=c\bar{s_2}^Yd, \ c,d\in Y^*$ and
$\delta(u\bar {s_1}^X), \ \delta(v\bar{s_2}^X)$ are disjoint. Then
there are two compositions according to $w^X=\delta(u\bar
{s_1}^X)a\delta(v\bar{s_2}^X)$ and
$w^X=\delta(v\bar{s_2}^X)a\delta(u\bar {s_1}^X)$ for $a\in X^*$. We
only prove the first.
\begin{eqnarray*}
(f,g)_w &=& \delta(us_1)a\delta(v\bar{s_2}^X)-\delta(u\bar{s_1}^X)ac\delta(vs_2)d \\
&\equiv&\delta(us_1\gamma(a)v\bar{s_2}^X-u\bar{s_1}^X\gamma(a)v\gamma(c)s_2\gamma(d)) \\
&\equiv&\delta(u\gamma(a)v(s_1\bar{s_2}^X-\bar{s_1}^X\gamma(c)s_2\gamma(d))) \\
&\equiv&\delta(u\gamma(a)vC_1(s_1,s_2,\gamma(c))_{w'}) \\
&\equiv&0 \ \ \ \ \ \ \ mod(S',w),
 \end{eqnarray*}
where
$w'=\bar{s_1}^X\bar{s_2}^X\bar{s_1}^Y=\bar{s_1}^X\bar{s_2}^X\gamma(c)\bar{s_2}^Y\gamma(d)$
and
$\overline{u\gamma(a)vC_1(s_1,s_2,\gamma(c))_{w'}}<u\gamma(a)vw'=\gamma(w)$.

\emph{(1.3)} $X,Y$-inclusion

We may assume that $\bar{g}$ is a subword of $\bar{f}$, i.e.,
$w=\delta(u\bar{s_1})=ac\delta(v\bar{s_2})bd$, $a,b\in X^*$, $c,d\in
Y^*$. Then $u\bar{s_1}^X=\gamma(ab)v\bar{s_2}^X=mL$ for some
$m\in[X]$, $u\bar{s_1}^Y=\gamma(c)\bar{s_2}^Y\gamma(d)$.
\begin{eqnarray*}
(f,g)_w &=& \delta(us_1)-ac\delta(vs_2)bd \\
&\equiv&\delta(us_1-\gamma(ac)vs_2\gamma(bd)) \\
&\equiv&\delta(m\frac{L}{\bar{s_1}^X}s_1-m\frac{L}{\bar{s_2}^X}\gamma(c)s_2\gamma(d)) \\
&\equiv&\delta(mC_1(s_1,s_2,\gamma(c))_{w'}) \\
&\equiv&0 \ \ \ \ \ \ \ mod(S',w),
 \end{eqnarray*}
where $w'=L\gamma(c)\bar{s_2}^Y\gamma(d)$ and
$\overline{mC_1(s_1,s_2,c)_{w'}}<mw'=\gamma(w)$.

$(2)$  $f \wedge h_{ij}$

Similar to the proof of Theorem \ref{th9}, they only have
compositions of $X$-intersection. Suppose that
$\bar{s_1}^X=x_{i_1}\cdots x_{i_r}x_i, \ (i_1\leq i_2\leq\cdots\leq
i_r\leq i)$. Then
$\bar{f}=\delta(u\bar{s_1})=x_{i_1}vx_i\bar{s_1}^Y$ for some $v\in
k\langle X\rangle, \ and \ v=\delta\gamma(v)$ and
$w=\delta(u\bar{s_1})x_j\bar{s_1}^Y.$

If $j\leq i_1$, then
\begin{eqnarray*}
(f,h_{ij})_w &=& \delta(us_1)x_j- x_{i_1}v\bar{s_1}^Y(x_ix_j-x_jx_i)\\
&=&\delta(u(s_1-\bar{s_1}))x_j+ x_{i_1}vx_jx_i\bar{s_1}^Y \\
&\equiv& x_j\delta(u(s_1-\bar{s_1}))+x_j x_{i_1}vx_i\bar{s_1}^Y \\
&\equiv& x_j(\delta(u(s_1-\bar{s_1}))+\delta(u\bar{s_1}) ) \\
&\equiv&x_j\delta(us_1)\\
&\equiv&0 \ \ \ \ \ \ \ mod(S',w).
 \end{eqnarray*}

If $j> i_1$, then $ux_j\in U(\bar{s_1})$ and
\begin{eqnarray*}
(f,h_{ij})_w &=& \delta(us_1)x_j- x_{i_1}v\bar{s_1}^Y(x_ix_j-x_jx_i)\\
&=&\delta(u(s_1-\bar{s_1}))x_j+ x_{i_1}vx_jx_i\bar{s_1}^Y \\
&\equiv& \delta(ux_j(s_1-\bar{s_1}))+\delta( x_{i_1}vx_ix_j\bar{s_1}^Y) \\
&\equiv& \delta(ux_j(s_1-\bar{s_1}))+\delta(ux_j\bar{s_1}) \\
&\equiv&\delta(ux_js_1)\\
&\equiv&0 \ \ \ \ \ \ \ mod(S',w).
 \end{eqnarray*}
This completes the proof.\ \ $\square$

\ \

\end{document}